\documentclass[mdpi,                
                preprintnumbers,
                amsmath,
                amssymb,
                showpacs,
                superscriptaddress]{revtex4-2}
\usepackage{graphicx} 
\usepackage{dcolumn}  
\usepackage{bm}       
\usepackage{amsfonts}
\usepackage{amssymb}
\usepackage[dvipsnames]{xcolor}
\usepackage{soul}
\usepackage{tabularx} 
\usepackage{booktabs}
\usepackage{braket}
\usepackage{amsmath}
\usepackage[normalem]{ulem}

\begin{document}

\title{EWS time delay in low energy $e-C_{60}$ elastic scattering}

\author{Aiswarya R.}
\affiliation{%
Department of Physics, Indian Institute of Technology Patna, Bihar 801103, India}

\author{Rasheed Shaik}
\affiliation{%
School of Physical Sciences, Indian Institute of Technology Mandi, Kamand, H.P. 175075, India}

\author{Jobin Jose}
\email[]{jobin.jose@iitp.ac.in}
\affiliation{%
Department of Physics, Indian Institute of Technology Patna, Bihar 801103, India}

\author{Hari R. Varma}
\affiliation{%
School of Physical Sciences, Indian Institute of Technology Mandi, Kamand, H.P. 175075, India}

\author{Himadri S. Chakraborty}
\affiliation{Department of Natural Sciences, D.L.\ Hubbard Center for Innovation, Northwest Missouri State University, Maryville, Missouri 64468, USA}

\maketitle
Time delay in a projectile-target scattering is a fundamental tool in understanding their interactions by probing the temporal domain. The present study focuses on computing and analyzing the Eisenbud-Wigner-Smith (EWS) time delay in low energy elastic $e-C_{60}$ scattering. The investigation is carried out in the framework of a non-relativistic partial wave analysis (PWA) technique. The projectile-target interaction is described in (1) Density Functional Theory (DFT) and (2) Annular Square Well (ASW) static model, and their final results are compared in details. The impact of polarization on resonant and non-resonant time delay is also investigated.
\section{Introduction}

The unprecedented developments in the field of attosecond chronoscopy have enabled scientists to follow the electron dynamics to a very high precision. This was possible by the pioneering experimental studies carried out by Pierre Agostini, Ferenc Krausz, and Anne L’Huillier for which they have been awarded the 2023 Nobel prize in Physics. The related comprehensive developments in research marked the birth of a new sub-field, "atto-science", in the vast landscape of ultrafast science. As a result, many of the fundamental interactions, initially considered as instantaneous processes, are found to be associated with some delay or acceleration of the order of attoseconds.  For example,  Krausz and coworkers \cite{schultze2010delay} have shown that the electron emission from the $3p$ orbital of the Ne atom has a temporal delay of $21\pm5$ attoseconds when compared to the delay from the $3s$ orbital. This finding instigated numerous theoretical and experimental investigations in time delay studies that explored features associated with photoionization dynamics.

Surprisingly not many studies, in the wake of these developments, were reported on the particle scattering processes in the temporal domain, with the exception of a few studies by Amusia \cite{amusia2018time_C60_sc, amusia2019partial, amusia2020time_endo}. The idea of time delay was initially proposed in the context of electron-scattering by Eisenbud \cite{eisenbud1948formal}, Wigner \cite{wigner1955lower}, and Smith \cite{smith1960lifetime}. Later a few notable studies were reported analyzing temporal dynamics in the scattering \cite{amusia2022wigner, amusia2018time_C60_sc,jose2022investigation}.  Eisenbud–Wigner–Smith (EWS) time delay is defined as the energy derivative of the scattering phase shift and it is an observable quantum mechanical parameter \cite{deshmukh2021time, schultze2010delay, guenot2012photoemission, wang2010attosecond}. The understanding of EWS time delay in a scattering process provides important details. For example, shape resonances in the scattering cross-section demonstrate how the projectile-target complex interacts to produce a quasi-bound state, which results in a temporal delay \cite{bain1972shape}. Studies on shape resonance have gained a lot of attention because of their widespread application in several fields of physics, including cold atom physics \cite{yao2019degenerate}, biological research \cite{bianconi1985co,bianconi1978k}, condensed matter physics \cite{bianconi1994possibility},  quantum transportation \cite{miroshnichenko2010fano}, etc. 

Over the years, many atomic and molecular targets have been analyzed using scattering techniques. Among these, fullerene- $C_{60}$ has been the subject of many spectroscopic investigations. Fullerene-$C_{60}$, having 12 pentagonal and 20 hexagonal carbon rings, is an extremely stable molecule with a fascinating symmetry. Another similar class of molecular compounds, endohedral fullerenes \cite{chai1991fullerenes}, are fullerene cages with an atoms or a molecule or a smaller fullerene inside them. These are denoted by $A@C_N$ where \textit{N} denotes the number of carbon atoms in the fullerene and \textit{A} is the species trapped in the cage. Fullerene and endofullerene complexes have a wide range of practical applications, such as in cancer detection and treatment \cite{hartman2008detecting}, medical imaging \cite{shu2008conjugation}, organic photovoltaic devices \cite{ross2009endohedral, hedley2013determining}, quantum computing \cite{harneit2007room}, for Hydrogen storage \cite{zhao2005hydrogen}, etc. Rubidium atom inside a fullerene, commonly known as Jahn-Teller metals, exhibits superconductivity at high-temperature \cite{zadik2015optimized}. In interstellar environments, traces of fullerenes and their ionic complexes have also been detected \cite{campbell2015laboratory, cami2010detection, cordiner2019confirming}.  

Despite such widespread applications, the majority of the work on fullerenes are photoionization-based (see articles \cite{pazourek2015attosecond, deshmukh2021time} and references therein). Recent studies have also addressed time delays in the photoionization of fullerene \cite{magrakvelidze2015attosecond, biswas2021attosecond} and endofullerene systems \cite{amusia2020time_ph_endo, dixit2013time_ph_endo}. However, in general, there are relatively limited theoretical and experimental studies investigating the elastic scattering off a fullerene. A few notable molecular-level calculations on $e-C_{60}$ elastic scattering were accomplished by the group of McKoy using Schwinger multichannel (SMC) method \cite{winstead2006elastic, hargreaves2010elastic} and by Gianturco et al. employing static exchange-correlation polarization (SECP) potential model \cite{gianturco1999computed}. On the experimental side, Tanaka et al. performed low-energy elastic $e-C_{60}$ scattering \cite{tanaka1994crossed, tanaka2021elastic}, whereas Hargreaves et al. investigated high-energy elastic scattering properties \cite{hargreaves2010elastic}. All of these aforementioned studies focused on total, partial, and differential cross-sections of scattering. No experimental elastic scattering time delay study of $C_{60}$ has been reported so far. On the theoretical side, a study by the group of Amusia et al. \cite{amusia2018time_C60_sc} used the Dirac bubble potential model to analyze the fullerene time delay, taking into account 6 partial waves. Additionally, a study conducted by Aiswarya and  Jose, employing only the annular square well (ASW) model potential, has analyzed EWS and angular time delay in $e-C_{60}$ scattering \cite{jose2022investigation}. 

The current work focuses on addressing the average time delay properties and the effects of polarization on time delay, using a more sophisticated modeling of fullerene based on Density Functional Theory (DFT) \cite{madjet2008photoionization, choi2017effects}. Additionally, we use the static ASW model potential for comparison \cite{dolmatov2014electron}. In order to investigate the resonant time delay features of the $e-C_{60}$ scattering, the low impact energy range is considered. In this energy range the projectile electron wavelength is so large that the electron diffraction effects can be neglected. Theoretical aspects of the work are covered in Section \ref{Sec: Theoretical Details}, the results and associated discussions are presented in Section \ref{Sec: Results}, and Section \ref{conculusion} concludes the report.

\section{Theoretical Details}
\label{Sec: Theoretical Details}
To simulate the $e-C_{60}$ interaction potential in the current investigation, two model potentials are used: (1) ASW model \cite{dolmatov2014electron} and (2) DFT potential generated within the Local Density Approximation (LDA) \cite{madjet2008photoionization, choi2017effects}. The ASW potential is expressed as \cite{dolmatov2014electron}:  
\begin{equation}
    V_{ASW}=\begin{cases}
         -U, & r_c-\frac{\Delta}{2}\le r\le r_c+\frac{\Delta}{2} , \\
          0, & otherwise
            \end{cases}
\end{equation}

where the cage thickness $\Delta$= 2.91 a.u., the mean radius \textit{$r_c$} and well depth are respectively 6.71 a.u. and 0.2599 a.u. The ASW parameters were selected to resemble some physical properties of fullerene closely: $\Delta$ is twice the covalent radius of carbon atom and well depth is adjusted to correspond to the electron affinity of -2.65 eV for $\ell=1$ state of $C_{60}$ \cite{dolmatov2014electron}. In effect, therefore, the ASW potential is static, being oblivious to C$_{60}$ electrons.  

\begin{figure}
\includegraphics[width=12.5 cm]{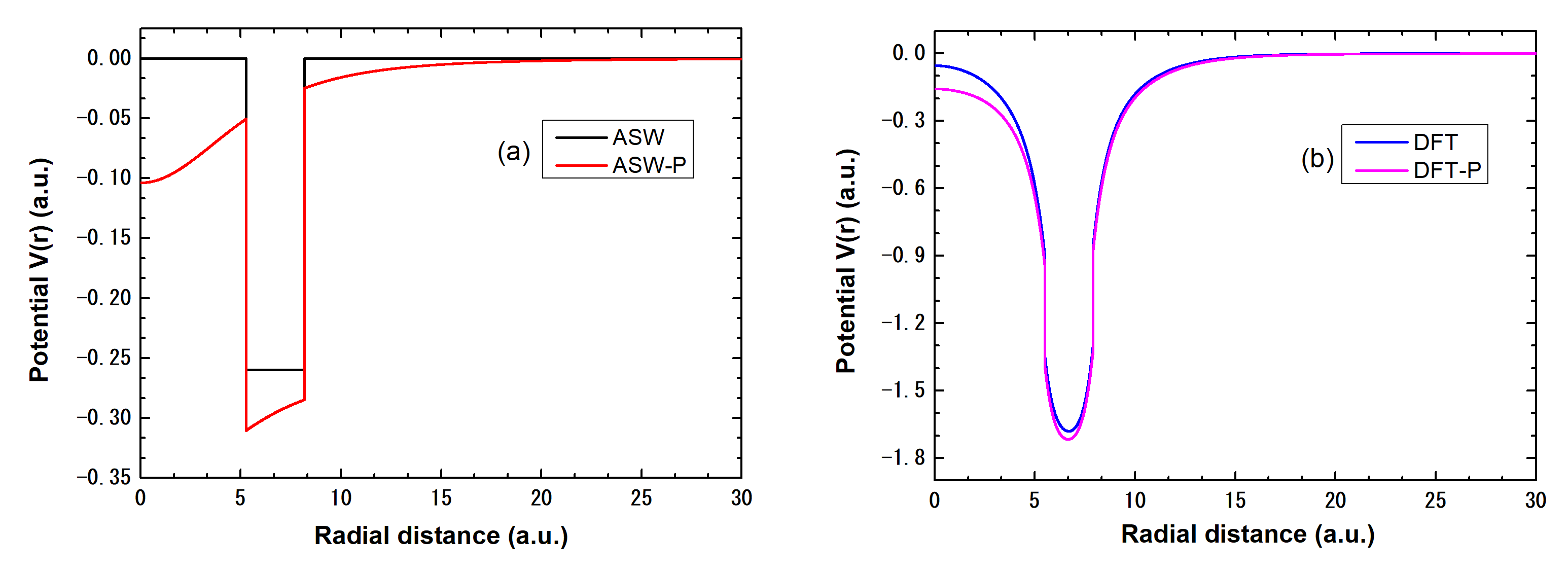}
\caption{ (a) ASW and (b) DFT model potentials with and without polarization effects.\label{fig1}}
\end{figure}   

In the DFT model, the complex core of fullerene is calculated within the local density approximation (LDA) model where the potential ($V_{jel}$) \cite{puska1993photoabsorption}, obtained after considering sixty $ C_4^+$ ions to be smeared and homogeneously distributed over a jellium shell of an average radius \textit{R} of 6.699 a.u. A constant pseudo-potential was added for the quantitative accuracy \cite{puska1993photoabsorption}. The fullerene shell then optimizes self-consistently (see below) to a thickness $\Delta$ of 2.411 a.u. and a value of the pseudo-potential to ensure the charge neutrality and to reproduce the measured first ionization threshold. The DFT-LDA potential for $e-C_{60}$ interaction can be written as \cite{madjet2008photoionization, choi2017effects}: 
\begin{equation}
\label{DFT_potential}
    V_{DFT}(r)=V_{jel}(r)+\int\frac{\rho(r^{\prime})}{r-r^{\prime}}dr^{\prime}+V_{xc}(\rho(r)),
\end{equation}
where the second term on the right side represents the direct interaction and the Kohn-Sham equations \cite{kohn1965self} for a system of 240 electrons are solved self-consistently to obtain the ground-state electron density $\rho(r)$ of the $C_{60}$ \cite{madjet2008photoionization, choi2017effects, madjet2010photoionization}. The third term in equation (\ref{DFT_potential}) corresponds to the exchange-correlation interaction. The specific form of exchange-correlation potential used in the present study is \cite{choi2017effects}:
\begin{equation}
  V_{xc}(\rho(r))=-\left(\frac{3\rho(r)}{\pi}\right)^{1/3}-0.0333\ln\left[1+11.4\left (\frac{4\pi\rho(r)}{3}\right)^{1/3}\right].
\end{equation}
Here, the first and second terms on the right-hand-side correspond to exchange and correlation potentials respectively. The exchange term is derived from the Hartree-Fock (HF) formalism \cite{choi2017effects} in LDA. Obviously, the DFT potential extends beyond the static model to incorporate the C$_{60}$ valence electron density within a mean field.

A realistic evaluation of the $e-C_{60}$ interaction requires accounting for polarization effects. As in the study by Dolmatov et al., the current study also uses a static dipole polarization potential of the form \cite{dolmatov2017effects, amusia2019angle_pol_exp_comp}:
\begin{equation}
\label{pol_pot}
    V_{pol}=\frac{-\alpha C_{60}}{2\left( r^2+b^2\right)^2}. 
\end{equation}
In equation (\ref{pol_pot}), the static dipole polarizability $\alpha =$850 a.u. and cut-off parameter \textit{b = }8 a.u., which approximately includes the fullerene's radial extent. The final effective potential after the addition of polarization potential is:
\begin{equation}
\label{eff_pot}
    V_{C_{60}}^{eff}=V_{DFT/ASW}+V_{pol}.
\end{equation}

To illustrate the effect of polarization, $V_{pol}$ is selectively included and omitted in the present study. Thus 4 case studies are carried out: (1) $ASW$ (ASW without polarization), (2) $ASW-P$ (ASW with polarization), (3) $DFT$ (DFT without polarization), and (4) $DFT-P$ (DFT with polarization). All the respective potential forms are shown in figure \ref{fig1}. The ASW potential form (figure \ref{fig1}(a)) is altered by the introduction of polarization potential, especially close to the fullerene center. The form of the DFT potential (figure \ref{fig1}(b)) is innately diffusive, and the addition of the attractive polarization potential further decreases the well depth. 

A non-relativistic partial wave analysis (PWA) is carried out for all the cases considered. The scattered wave function of the electron after interacting with the $C_{60}$ potential will be the solution of the following radial Schr\"{o}dinger equation: 
   
\begin{equation}
\label{sch_eq}
    \left\{\ \frac{-\hbar^2}{2m}\frac{d^2}{dr^2}+\left[ V_{C_{60}}^{eff}+\frac{\hbar^2\ell(\ell+1)}{2mr^2}\right] \right\}\ u_\ell(r)=E u_\ell(r).
\end{equation}
Separate calculations are performed after decoupling the effect of polarization to investigate its impact. Using Numerov's method \cite{thijssen2007computational}, the solution of equation (\ref{sch_eq}) is computed with the proper boundary conditions. Care is taken to ensure the continuity of the wave function and its derivative at the boundary of $C_{60}$. As in the earlier study \cite{jose2022investigation}, the elastic scattering phase shift ($\delta_\ell$) of the $\ell^{th}$ partial wave is obtained by considering the asymptotic form of the wave function \cite{joachain1975quantum}: 
\begin{equation}
\label{asymp_bhr}
    u_\ell(r>r_{max})\propto kr[j_\ell(kr) cos\delta_\ell -n_\ell(kr) sin\delta_\ell ].
\end{equation}

For any energy \textit{E} of the projectile, the wave vector $k=\frac{\sqrt{2mE}}{\hbar}$. In equation (\ref{asymp_bhr}), $j_\ell$ and $n_\ell$ represent the Bessel function of the first and second kind, respectively, and $r_{max}$ = 28 a.u. is taken as the practical infinity. Let $r_1$ and $r_2$ be the two radial points beyond $r_{max}$, using which phase shift ($\delta_\ell$) can be computed as \cite{joachain1975quantum}:

\begin{equation}
    tan \delta_\ell(k)=\frac{\zeta j_\ell(kr_1)-j_\ell(kr_2)}{\zeta n_\ell(kr_1)-n_\ell(kr_2)}
\end{equation}
with
\begin{equation}
    \zeta=\frac{r_1 u_\ell(r_2)}{r_2 u_\ell(r_1)}.
\end{equation}
Numerical values of $j_\ell$ and $n_\ell$ are generated using well-established subroutines SPHJ and SPHY \cite{zhang1997computation}. Partial waves with $\ell=$0–15 were determined to be sufficient for both model potentials up to the energy of \textit{E = }0.5 a.u. considered in this study. Since most of the prominent resonances are observed in this energy range, the low energy elastic $e-C_{60}$ scattering is addressed. The total cross-section (TCS) of elastic scattering is given by \cite{thijssen2007computational, joachain1975quantum}: 

\begin{equation}
    \sigma_{total}=\frac{4\pi}{k^2} \sum_{\ell=0}^{\infty} (2\ell+1) sin^2\delta_\ell.
\end{equation}
Analysis of the quasi-bound states arising from distinct resonant partial-wave interactions with fullerene is performed using the Fano parametrization \cite{fano1961effects} formula: 

\begin{equation}
\label{fano}
    \sigma_r=\frac{\sigma_o (q+\varepsilon)^2}{1+\varepsilon^2},
\end{equation}
with
\begin{equation}
\label{varepsilom}
    \varepsilon=\frac{E-E_r}{\Gamma/2}.
\end{equation}
In equation (\ref{fano}), \textit{q} is the shape parameter, $\sigma_r$ and $\sigma_o$ are the resonant and background cross sections respectively. Equation (\ref{varepsilom}) provides the parameter $\varepsilon$, where \textit{E} is the incident projectile energy and $\Gamma$ and $E_r$ are the resonant width and energy, respectively. The time delay $\tau_\ell$ of the $\ell^{th}$ partial wave is given by the EWS time delay formula \cite{smith1960lifetime}:  
\begin{equation}
\label{time_delay}
    \tau_\ell(E)=2\hbar\frac{\partial \delta_\ell}{\partial E}=2\hbar\delta'_\ell.
\end{equation}
The average time delay for a given potential is provided by \cite{amusia2018time_C60_sc}:
\begin{equation}
\label{avg_td}
    \tau_{avg}(E)=\sum_{\ell=0}^{\infty}\frac{\sigma_\ell}{\sigma_{total}}\tau_\ell(E).
\end{equation}
where $\sigma_{\ell}$ and $\tau_\ell$ are respectively, the partial cross-section (PCS) and time delay of the $\ell^{th}$ partial wave. We analyze the resonant time delay for all 4 case studies of model potentials (see above). All the calculations are done in atomic units (a.u.) unless specified otherwise.
\section{Results and Discussion}
\label{Sec: Results}
 The section is divided into 4 subsections: section (\ref{Sec: TCS, low energy behavior}) deals with the $e-C_{60}$ TCS and near-threshold behavior of phase shift and time delay, and section (\ref{Sec: Fano parameterization of resonance}) discusses the Fano parameterization of resonant partial waves. Resonant phase shifts and time delays are covered in detail in section (\ref{Sec: Resonant phase shift and time delay}), and section (\ref{Sec: Average time delay}) presents the average time delay. 

\subsection{TCS and near-threshold behavior of phase shift and time delay}
\label{Sec: TCS, low energy behavior}

\begin{figure}
\includegraphics[width=12.5 cm]{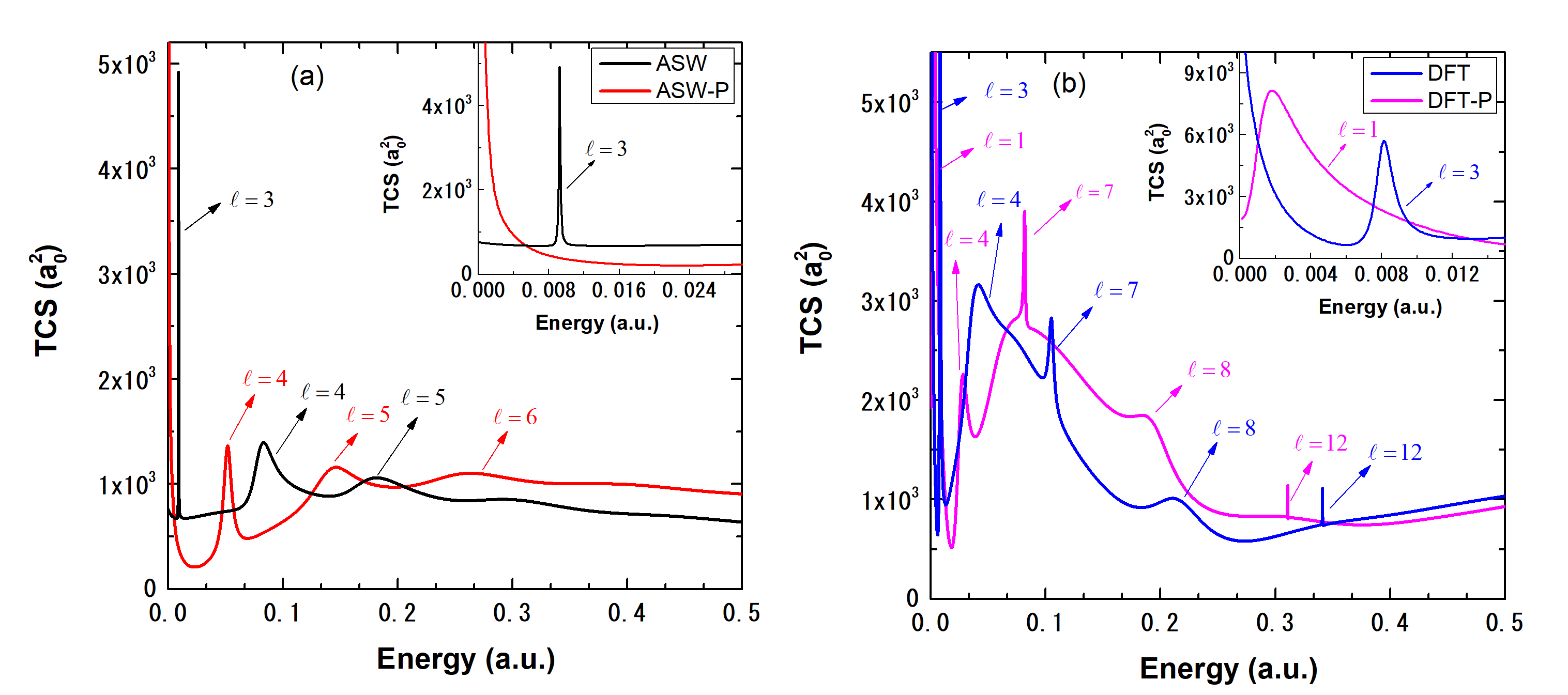}
\caption{ TCS of $e-C_{60}$ elastic scattering using (a) ASW (black) and ASW-P (red); (b) DFT (blue) and DFT-P (magenta). Resonances are labeled with corresponding partial wave numbers. The inset shows magnified TCS in the low-energy region.\label{fig2}}
\end{figure}

The TCS of $e-C_{60}$ elastic scattering contributed by $\ell$ = 0-15 partial waves is shown in figure \ref{fig2}. Cross-sections in the ASW  and ASW-P model are shown in figure \ref{fig2}(a).  In TCS, for the ASW case, 3 prominent peaks are noted which are from the partial waves $\ell$=3, 4, and 5; $\ell=3$ is the sharpest and closer to zero energy, while $\ell=5$ is the weakest. This sharp resonance corresponding to $\ell=3$ is not observed in the case of polarization-added ASW potential case (ASW-P). In the ASW-P case, resonant peaks are mainly due to the contribution from partial waves $\ell$= 4, 5, and 6. The inset features a magnified view of TCS at the near-zero energy region which shows the polarization induced  $\ell=0$ resonance, that occurs right at the threshold. For ASW potential, the addition of polarization effects shifts the resonance to lower energies. Thus, two common resonance peaks ($\ell=$ 4 and 5) are noticed both in ASW and ASW-P. The TCS obtained using DFT and DFT-P potential models are shown in figure \ref{fig2}(b). Since the DFT simulated $e-C_{60}$ interaction potential has a deeper well depth than the ASW case, more resonances are formed. A zero-energy resonance is seen for the DFT potential, followed by a resonance corresponding to $\ell= 3$. For the polarization added DFT model (DFT-P), a resonance for $\ell = 1$ is noted quite close to the threshold. A set of peaks corresponding to the partial waves $\ell=$ 4, 7, 8, and 12 is seen both in DFT and DFT-P. When the polarization is added to the model potentials (ASW/DFT), therefore, a consistent trend is that the resonant partial wave peaks shift towards lower energies. This is because the presence of the attractive polarization potential results in a minor reduction in the angular momentum barrier, leading to the attainment of resonant conditions at slightly lowered projectile energies. 

\begin{figure}
\centering
\includegraphics[width=18cm]{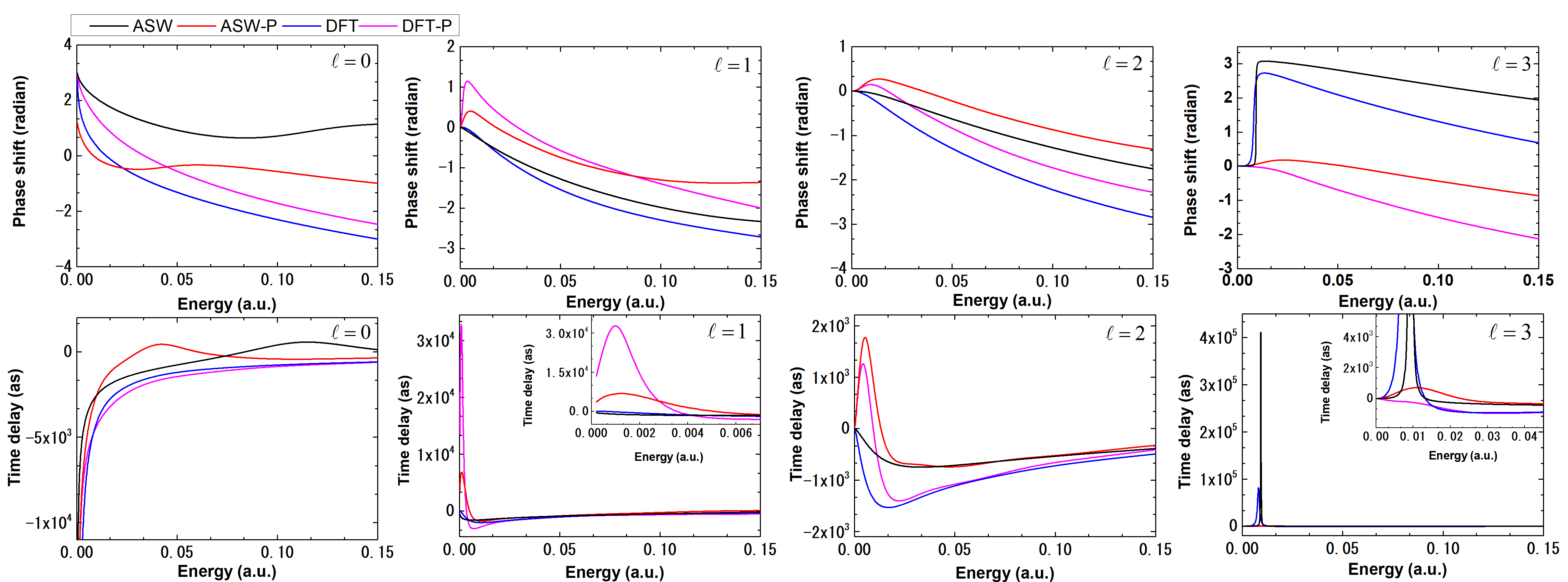}
\caption{Phase shift (upper panel) and corresponding time delay (lower panel) in the near-zero energy range using model potentials ASW (black) ASW-P(red), DFT (blue) and DFT-P (magenta). Inset shows the magnified view of time delay for $\ell$=1 and 3 partial waves. \label{fig3}}
\end{figure}  

Figure \ref{fig3} illustrates the near-threshold phase shift and EWS time delay behavior for ASW and DFT approximations; time delays are presented in attoseconds (as). In the zero energy limit, the phase shift and time delay follow the Wigner threshold law \cite{wigner1948behavior_threshold_law}. Since the \textit{s-}wave phase shift behaves as $\delta_{\ell=0}(E\to 0)\propto \pi-E^{1/2}$, it tends to $\pi$ as $E\rightarrow0$, which is evident in the figure for the case of $\ell=0$. The corresponding time delay goes to negative infinity $\tau_{\ell=0}(E\to 0)\propto -E^{-1/2}$, which is also clear in the bottom panel for $\ell=0$. According to the Wigner threshold law, the phase shift for the other partial waves $\ell \ne 0$ goes as $\delta_{\ell>0}(E\to 0)\propto \pm E^{\ell+1/2}$. Therefore, the phase shift tends to zero in the zero-energy limit, as shown in the upper panel in figure \ref{fig3}. The corresponding  time delay vanishes for the partial waves $\ell>0$ in the low energy limit, since $\tau_{\ell>0}(E\to 0)\propto \pm E^{\ell-1/2}$.  Hence, the current set of results is consistent with the Wigner threshold law. 

We now investigate the low-energy behavior of phase shift and time delay upon inclusion of the polarization effect. For each of $\ell=$ 1 and 2, the phase shift is positive in the low-energy region in ASW-P and DFT-P compared to the non-polarized target approximation. Consequently, there is a positive hump in the time delay for both these cases. We note that the polarization induces an attractive potential (figure \ref{fig1}). Let us suppose that the polarization depends parametrically on a quantity $\lambda$, which can be tuned to adjust the polarization potential from the unpolarized state. In the limit $\delta_\ell \to 0$, the relation between the phase shift and the potential can be approximated as \cite{joachain1975quantum}:
\begin{equation}
\label{ps_V_relation}
    \frac{d\delta_\ell}{d\lambda}=-k\int_{0}^{\infty} [u_\ell(\lambda,r)]^2 \frac{\partial V_{C_{60}}}{\partial \lambda} \,dr,
\end{equation}
where $\lambda$ is the tuning parameter mentioned above and $u_\ell(\lambda,r)$ is the asymptotic radial wave function for wave vector \textit{k}. The rate of change of the strength of the $C_{60}$ potential upon the addition of the polarization effect is expressed by $\frac{\partial V_{C_{60}}}{\partial \lambda}$, which is negative in the present case. From equation (\ref{ps_V_relation}), it is evident that the associated phase shift variation has the opposite sign if the potential variation, $\frac{\partial V_{C_{60}}}{\partial \lambda}$, has the same sign for all values of \textit{r} \cite{joachain1975quantum}. Analyzing the phase shift behavior for $\ell=$1 and $\ell=$2, we find that phase shifts from all 4 model potentials approach zero in the very low energy range, validating the Wigner threshold law. Furthermore, for model potentials where polarization is considered (ASW-P/ DFT-P), the phase shift develops in positive values, and for the bare cases (ASW/  DFT), it grows in negative values.  The reason for this is found clear from equation (\ref{ps_V_relation}). Analyzing the polarized target situation (P case),  the value of $\frac{\partial V_{C_{60}}}{\partial \lambda}$ is negative compared with the non-polarized case; hence, at a given \textit{k} value, the scattering phase shift of the P case will be relatively positive compared to the non-P case. Coming to the behavior of the time delay, we recall that $\tau_\ell$ is obtained as the energy derivative of the $\delta_\ell$. As a result, for $\ell=$1 and $\ell=$2, the time delay is negative for the non-P case and positive for the P case. As $k$ increases, however, the equation (\ref{ps_V_relation}) ceases to be valid and the evolution of the phase shift begins to depend numerically on the strength and relative difference of the potentials. A peak in time delay is noticed as a result of the inclusion of polarization in the low-energy region. For the partial wave $\ell=3$, a completely different behavior of phase shit and time delay is obtained: The resonance condition is satisfied for both the ASW and DFT model potentials. In addition, an abrupt phase leap through $\pi/2$ and an overall shift by $\pi$ are noticed. Thus, a sharp peak in time delay is found. However, for ASW-P and DFT-P, the corresponding resonances have moved below the threshold as ascertained earlier, and hence no abrupt change in phase shift is observed. As a result, the time delay profile has no sharp resonance peaks. A detailed analysis of the resonant cross-section is presented in section \ref{Sec: Resonant phase shift and time delay}. As the $\ell$ value increases, the angular momentum barrier increases, which prohibits the intrusion of scattered waves in the region close to the origin. Consequently, the difference in the potential for the P and non-P cases will be less conspicuous. 

\subsection{Fano parameterization of resonance}
\label{Sec: Fano parameterization of resonance}

Fano parametrization is employed for all model potentials to determine the exact resonant energy $E_r$ and other resonant parameters \cite{fano1961effects}. The Fano resonance formula provides a way to analyze resonant scattering and was introduced by Ugo Fano in 1961 \cite{fano1961effects}. The shape parameter $q$ (from equation (\ref{fano})) characterizes the asymmetry and shape of the resonance profile. It describes the interference between a pure resonant channel (with a Lorentzian shape profile) and a non-resonant background. Resonances are of different shapes based on the altering strength of this interference. 
  
The Fano shape parameters for the resonance in the ASW $(\ell=3,4,5)$ and ASW-P $(\ell=4,5,6)$ cross-section are shown in table \ref{table1}. The corresponding resonant partial wave cross-sections are compared with the fitted cross-section in figure \ref{fig4}. The $q$ parameter value is large ($q=$ 28) and positive for $\ell=3$ partial wave, and has a lower background cross-section ($\sigma_o$= 8.00 $a_0^2$) for the ASW potential. A high $q$ value suggests that the interference between resonant and direct (background) scattering is weak, leading to a Lorentzian shape of the cross-section.  Furthermore, the weak coupling between the background and the resonant channel is reflected in the low $\sigma_o$. Relative to the other resonances ($\ell=$ 4 and 5 cases), the $\ell=3$ resonance shape is narrower $(\frac{\Gamma}{2}=0.08 \times10^{-3}$ a.u.) and more akin to the Lorentzian profile in the ASW case, as seen in figure \ref{fig4}(upper panel). Additionally, a small resonance width for $\ell=$ 3 indicates a longer resonance lifetime. The $q$ value remains positive for $\ell=$ 4 and 5 cases also, however, the values are seen to approach zero. A lower value of $q$ indicates a strong interference between resonant and direct scattering channels leading to asymmetric resonant shape profiles, as clear from figures \ref{fig4}(upper panel; $\ell=$ 4 and 5). Furthermore, the $\ell=4$ and $\ell=5$ resonances are weak, indicative of a shorter lifetime. 
  
For ASW-P (figure \ref{fig4}(lower panel)), no resonance is obtained for $\ell=$ 3. From table \ref{table1}, it can be noted that when polarization is introduced to ASW potential, the $q$ value for $\ell=$ 4 partial wave decreases from 2.69 to -19.30. Also, a considerable lowering of $\sigma_o$ and resonance width values is noted. This indicates that the $\ell=4$ partial wave becomes sharper with the addition of polarization. For $\ell=$ 5, the $q$ value goes from 1.19 to 3.90 when the polarization is added and a decrease in $\sigma_o$ and $\Gamma$ values is also observed. The resonance for $\ell=$ 6, as emerged for ASW-P, has a $q$ value of 1.69. Also, among the resonances observed in the ASW-P cases, the maximum $\sigma_o$ and minimum $\Gamma$ values are observed for $\ell=$ 6, indicating an asymmetric and weaker resonance.
The Fano resonance profiling is done for the  DFT and DFT-P model potential cases as well; the resulting values of the Fano parameters are listed in table \ref{table2}. In figure \ref{fig5}, the fitted profile for DFT (upper panel) and DFT-P (lower panel) are shown. The fitting parameter values have evolved quite differently between DFT and DFT-P cases. In DFT, a resonance is noted for the partial wave $\ell=$ 3, which has a $q$ value of 6.23 and a considerably large background cross-section, $\sigma_o= 135 a_0^2$. Consequently, the resonant cross-section is asymmetric. A resonance exists for $\ell=$ 1 in DFT-P with resonant energy $E_r=$ 0.0013 a.u.,  $q=$ 2.48 and $\sigma_o= 920 a_0^2$, and its asymmetry is indicative of the larger background cross-section. A common set of resonant partial waves ($\ell=$ 4, 7, 8, and 12) is observed for DFT and DFT-P cases. The $q$ value increases from 1.65 to 2.79 for $\ell=$ 4 when the polarization effect is added. Furthermore, a decrease both in the resonant width and background cross section is seen, indicating that for $\ell=4$ the polarization increases the sharpness of the resonance. When polarization effects are taken into account, the $q$ value for $\ell=$ 7 rises from -21.20 to -9.30, while for $\ell=$ 8, the $q$ value falls from 12.45 to -15.25. From the data table \ref{table2}, it could be noted that the resonance width is the narrowest for $\ell=$ 12 in both DFT and DFT-P. A narrow resonance width is indicative of a prolonged resonance time, as was previously mentioned. The corresponding $q$ values for DFT and DFT-P are -5.60 and -4.53, respectively. 

We remark that the lifetimes of these resonances, as derivable from their line-widths, will affect the scattering time delay in general, since the electron will be transiently captured in these resonance states before scattering to the detector.

\begin{table} 
\begin{center}
\caption{Fano parameters of resonances using ASW and ASW-P model potential.\label{table1}}
\newcolumntype{C}{>{\centering\arraybackslash}X} 
\begin{tabularx}{\textwidth}{C C C C C C}
\toprule
 &  & \multicolumn{4}{c}{\textbf{ Fano parameters}}	\\

\textbf{Model potential} & $\ell$ & ${E_r}$ & \textit{q} & $\frac{\Gamma}{2}$ & $\sigma_o$ \\
 & & a.u. & & $\times10^{-3}$ a.u. &  $a_0^2$ \\                           
\cline{1-6}
ASW    & 3 & 0.0092     & 28.00      & 0.08   & 8.00 \\
       & 4 & 0.0795     & 2.69       & 10.40  & 82.00 \\
       & 5 & 0.1620     & 1.19       & 31.00  & 148.00 \\
       &   &            &            &        &     \\
ASW-P & 4 & 0.0524     & -19.30     & 4.40   & 2.90 \\
       & 5 & 0.1380     & 3.90       & 26.00  & 29.00 \\
       & 6 & 0.2350     & 1.69       & 60.00  & 76.00 \\
\cline{1-6} 
\end{tabularx}
\end{center}
\end{table}

\begin{figure}
\includegraphics[width=14.5 cm]{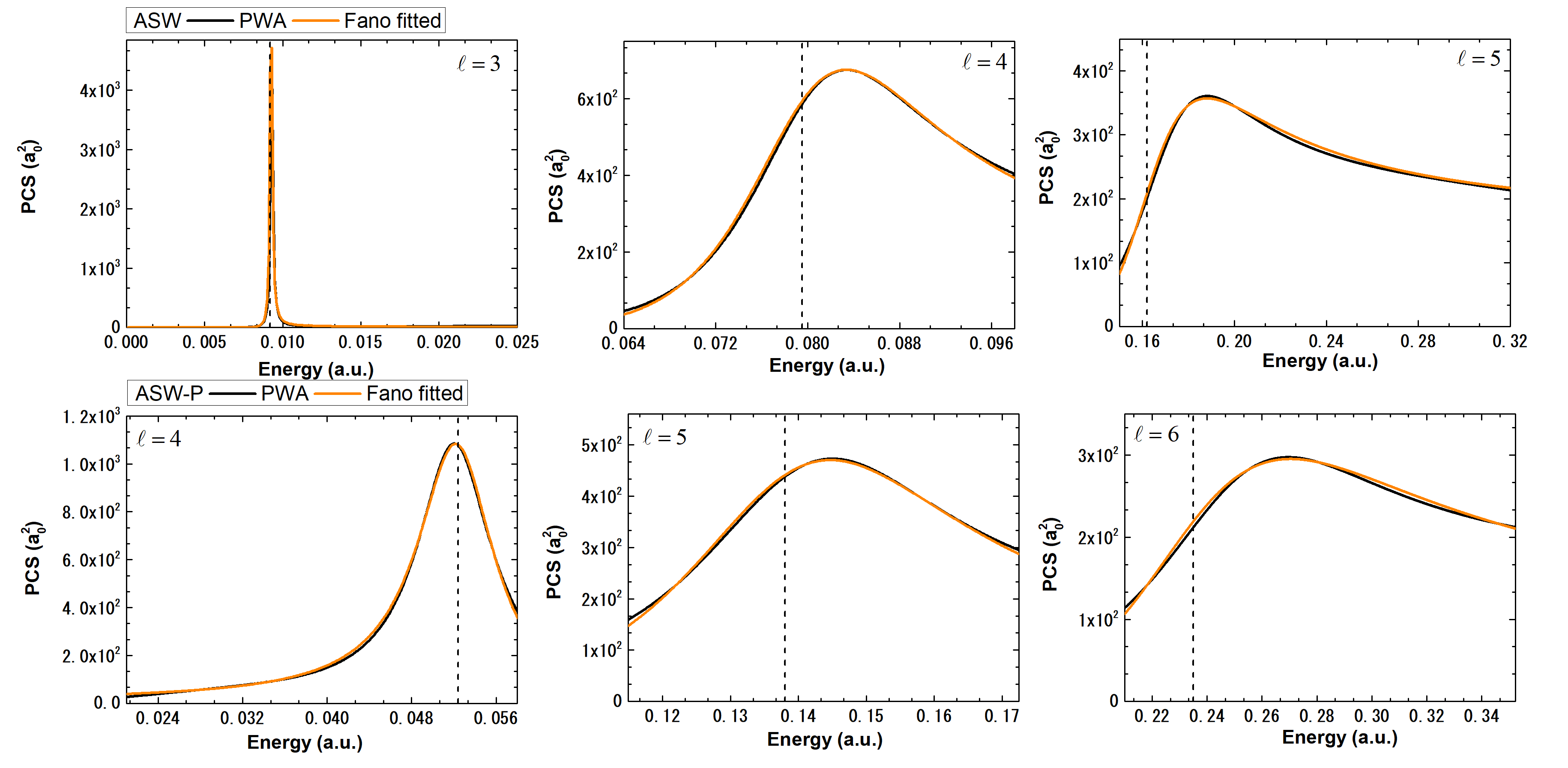}
\caption{ Comparison of $e-C_{60}$ elastic scattering resonant cross-sections calculated using PWA (black) with that from Fano fitting (orange) for ASW (upper panel) and  ASW-P (lower panel). Vertical line indicates the resonant energy $E_r$.\label{fig4}}
\end{figure}   

\begin{table}
\begin{center}
\caption{Fano parameters of resonances using DFT and DFT-P model potential.\label{table2}}
\newcolumntype{C}{>{\centering\arraybackslash}X} 
\begin{tabularx}{\textwidth}{C C C C C C}
\toprule
 &  & \multicolumn{4}{c}{\textbf{ Fano parameters}}	\\
\textbf{Model potential} & $\ell$ & ${E_r}$ & \textit{q} & $\frac{\Gamma}{2}$ & $\sigma_o$ \\
 & & a.u. & & $\times10^{-3}$ a.u. &  $a_0^2$ \\                       
\cline{1-6} 
DFT    & 3  & 0.0080     & 6.23    & 0.60   & 135.00 \\
       & 4  & 0.0340     & 1.65    & 11.40  & 359.00 \\
       & 7  & 0.1050     & -21.20  & 2.50   & 2.00 \\
       & 8  & 0.2120     & 12.45   & 28.20  & 3.20 \\
       & 12 & 0.3414     & -5.60   & 0.09   & 12.00 \\
       &    &            &         &        &         \\ 
DFT-P & 1  & 0.0013     & 2.48    & 1.35   & 920.00 \\
       & 4  & 0.0258     & 2.79    & 5.50   & 228.00\\
       & 7  & 0.0818     & -9.30   & 0.75   & 13.10 \\
       & 8  & 0.1920     & -15.25  & 20.60  & 2.40 \\
       & 12 & 0.3109     & -4.53   & 0.50   & 16.70 \\
\cline{1-6} 
\end{tabularx}

\end{center}
\end{table}

\begin{figure}
\centering
\includegraphics[width=18.5cm]{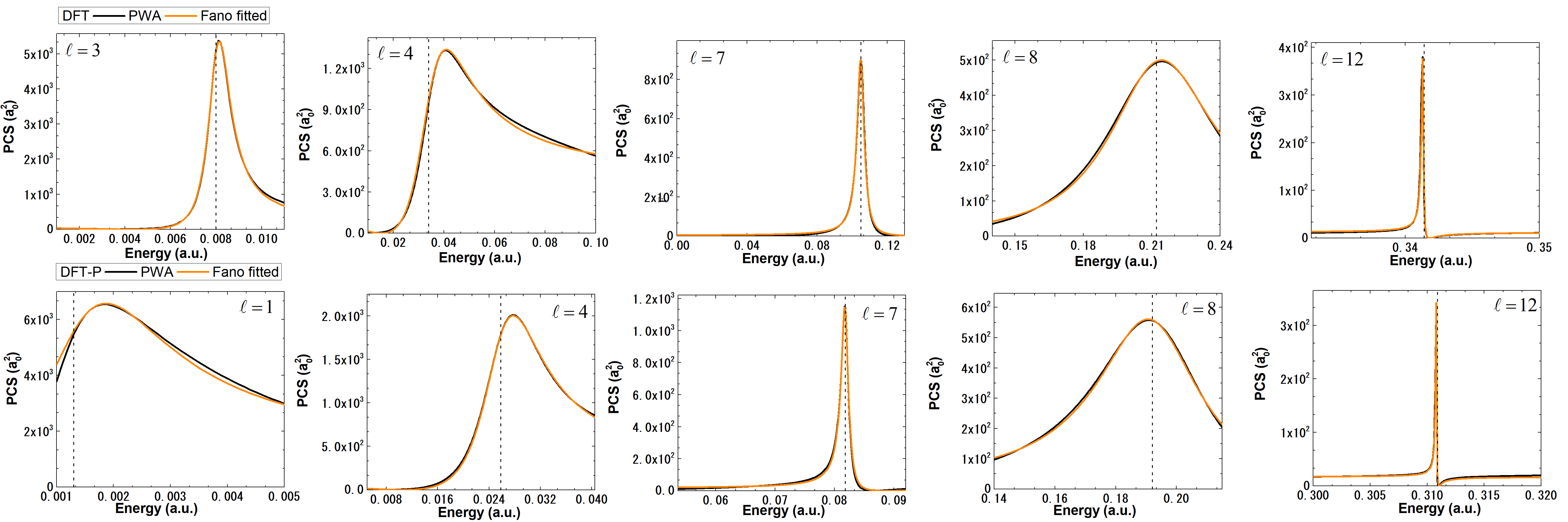}
\caption{Comparison of $e-C_{60}$ elastic scattering resonant cross-sections calculated using PWA (black) with that from Fano fitting (orange) for DFT (upper panel) and  DFT-P (lower panel). Vertical line indicates the resonant energy $E_r$.\label{fig5}}
\end{figure}

\subsection{Resonant phase shift and time delay}
\label{Sec: Resonant phase shift and time delay}

\begin{figure}
\includegraphics[width=14.5 cm]{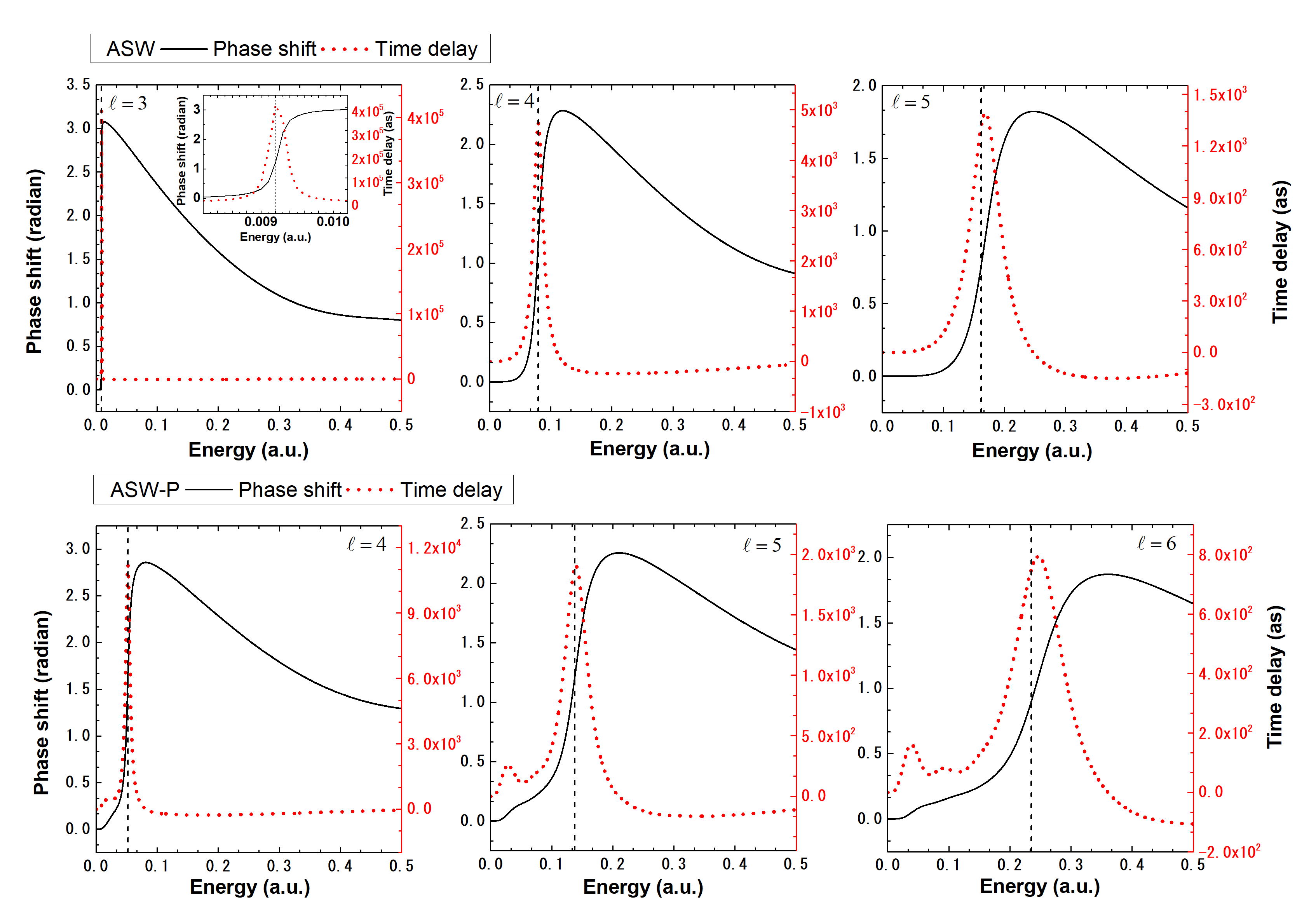}
\caption{ Resonant partial wave scattering phase shift (solid black) and EWS time delay (dashed red)  plotted for ASW (upper panel) and ASW-P (lower panel) potentials. The phase shift is linked to the left axis and the time delay is linked to the right axis.  Vertical line indicates the resonant energy $E_r$.\label{fig6}}
\end{figure}

In figure \ref{fig6}, the scattering phase shift and EWS time delay (using equation (\ref{time_delay})) are plotted for the resonant partial waves of ASW (upper panel) and ASW-P (lower panel) model potentials. The vertical lines indicate the corresponding resonant energies, which are obtained from the Fano analysis. First, we discuss the Non-P case (figure \ref{fig6}(upper panel)). The phase shift for $\ell=3$ partial wave shows an abrupt $\pi/2$ jump and a sweep by $\pi$ radian within a short energy range. As a result, the time delay sharply peaks as shown in the figure \ref{fig6}(upper panel; $\ell=3$). Subsequent partial waves ($\ell=$ 4 and 5) also meet the resonant condition but comparatively weakly. For these weaker resonances, the magnitude of the phase shift is less than $\pi$ across the resonance. The resonant time delay for the $\ell=3$ partial wave is 410456.2 as. This indicates that at this resonance the electron stays in the vicinity of the scatterer for quite a while. Resonant energy and the corresponding time delay values for the resonant partial waves are listed in table \ref{table3}. According to table \ref{table3}, the resonant time delay drops gradually for partial waves having higher angular momentum quantum numbers. A study by Amusia et al. \cite{amusia2018time_C60_sc}, using the Dirac bubble potential model to describe fullerene, reported a similar trend in the time delay profile. In the ASW-P case (figure \ref{fig6}(lower panel)), $\ell=3$ partial wave is non-resonant, but the partial waves $\ell=$4, 5, and 6 are resonant. The strongest resonance is observed for $\ell=4$, having a time delay value of 11092.1 at the resonant energy $E_r=$0.0524 a.u., followed by a smaller resonant time delay of 1893.9 as for $\ell=5$ and 751.1 as for $\ell=6$. A progressive reduction in the resonant time delay is observed for larger angular momentum partial waves, similar to the trend in ASW case. 

As analyzed in section \ref{Sec: TCS, low energy behavior}, a small hump appears in the phase shift near zero energy due to the polarization effects. Consequent to this, a smaller peak in the time delay is seen in the ASW-P case, which is absent in the ASW case. For $\ell=4$ partial wave, the peak in time delay due to polarization is at 0.0524 a.u., and for $\ell=5$ and $\ell=6$, respectively, at 0.1380 a.u. and at 0.2350 a.u.  Noticeably, although the phase shift is less sensitive to the polarization effect, an amplified feature is noted in time delay. This suggests that the time delay is more sensitive to changes in the potential compared to other scattering parameters.

\begin{table}
\begin{center}
\caption{Resonant energy and time delay for ASW and ASW-P case.\label{table3}}
\newcolumntype{C}{>{\centering\arraybackslash}X} 
\begin{tabularx}{\textwidth}{C C C C C}
\toprule
\textbf{Resonant partial wave}& \multicolumn{2}{c}{\textbf{ASW}}	& \multicolumn{2}{c}{\textbf{ASW-P}}\\
 $\ell$   & \textbf{Energy (a.u.) }& \textbf{Time delay (as)} & \textbf{Energy (a.u.)}& \textbf{Time delay (as)} \\
 \cline{1-5}                              
\midrule
3		& 0.0092   & 410456.2 & --      & --\\
4		& 0.0795   & 4731.3   & 0.0524	& 11092.1 \\
5		& 0.1620   & 1320.9   & 0.1380  & 1893.9 \\
6		& --       & 	--    & 0.2350	& 751.1 \\
 \cline{1-5}   
\end{tabularx}
\end{center}
\end{table}

\begin{figure}
\centering
\includegraphics[width=18cm]{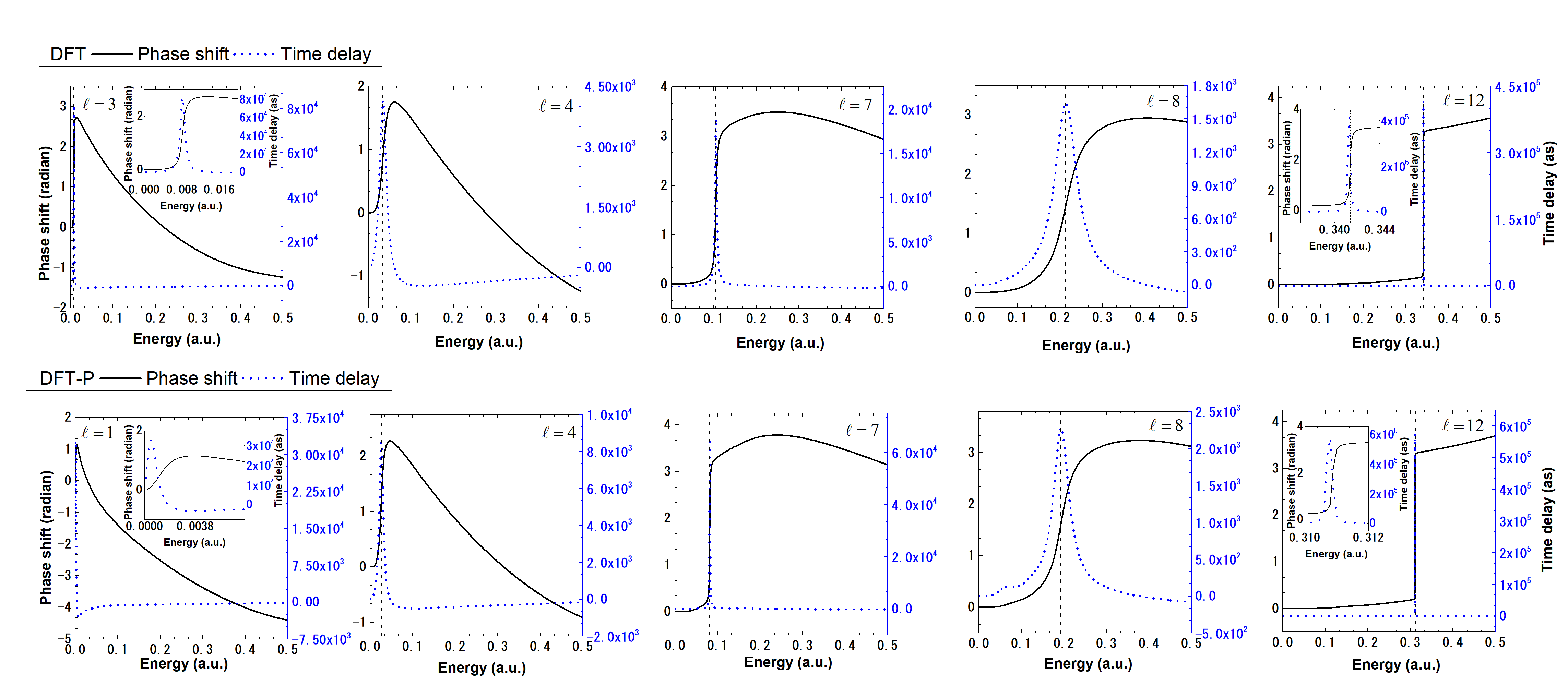}
\caption{Resonant partial wave scattering phase shift (Solid black) and EWS time delay (dashed blue)  plotted for DFT (upper panel) and DFT-P (lower panel) potentials. The phase shift is linked to the left axis and the time delay is linked to the right axis. Vertical line indicates the resonant energy $E_r$. \label{fig7}}
\end{figure} 

\begin{table} 
\begin{center}
\caption{Resonant energy and time delay for DFT and DFT-P case.\label{table4}}
\newcolumntype{C}{>{\centering\arraybackslash}X} 
\begin{tabularx}{\textwidth}{C C C C C}
\toprule
\textbf{Resonant partial wave}& \multicolumn{2}{c}{\textbf{DFT}}	& \multicolumn{2}{c}{\textbf{DFT-P}}\\
 $\ell$   & \textbf{Energy (a.u.) }& \textbf{Time delay (as)} & \textbf{Energy (a.u.)}& \textbf{Time delay (as)} \\
\cline{1-5}                                 
\midrule
1		& --     & --         & 0.0013  & 29683.7\\
3		& 0.0080 & 	81602.4   & --	    & -- \\
4		& 0.0340 & 	4128.5    & 0.0258  & 8493.2 \\
7		& 0.1050 & 	18799.3   & 0.0818	& 65485.2 \\
8       & 0.2120 &  1634.1    & 0.1920  & 2264.1 \\
12		& 0.3414 & 	227293.7 & 0.3109	& 237539.9 \\
\cline{1-5}   
\end{tabularx}
\end{center}
\end{table}

For the $e-C_{60}$ interaction modeled by DFT, more resonances are seen. This is, as noted before, because of the intrinsic diffusive shape and a deeper well in the DFT model. Phase shifts and corresponding time delay profiles of resonant partial waves are shown in figure \ref{fig7}. We first analyze the DFT case without polarization (refer figure \ref{fig7}(upper panel)). Here, the resonance condition is satisfied for the partial waves $\ell=$ 3, 4, 7, 8, and 12. Partial wave with $\ell=$ 12 exhibits the sharpest resonance, followed by partial waves with $\ell=$3, 7, 4, and 8. The longest resonant time delay is noted at $E_r=$ 0.3414 a.u. for $\ell=$ 12 with $\tau_{\ell=12}= 227293.7$ as. At $E_r=$ 0.2120 a.u., the shortest resonant time delay of 1634.1 as is observed for $\ell=8$ partial wave. Details of resonant energy and time delay values are given in table \ref{table4}. From figure \ref{fig7}(lower panel), the polarization added DFT model (DFT-P) exhibits 5 notable resonances, just as the DFT case. Here, a near-threshold resonance is obtained for $\ell=1$, for which the resonant energy and time delay values are 0.0013 a.u. and 29683.7 as respectively. Next resonance condition is attained for partial wave $\ell=4$, followed by $\ell=$ 7, 8, and 12. Similar to the DFT case, a steep jump in phase shift is obtained for $\ell=12$, producing the maximum resonant time delay value of 237539.9 as at the energy $E_r=$ 0.3109 a.u. A steady jump in the phase is seen for the partial wave with $\ell=7$ which induces a time delay value of 65485.2 as, while for $\ell=4$ the resonant time delay is 8493.2 as. In comparison, the phase shift change is relatively small for $\ell=8$ resulting in the lowest resonant time delay of 2264.1 as. For the polarization-added case, however, there is a slight decrease in the resonant energy compared to the non-P case. This is because, as previously mentioned, the attractive polarization potential lowers the barrier potential. Similar to the ASW analysis, a very small bump in time delay for DFT-P, for $\ell=8$ is found which is likely induced by the polarization effect of $C_{60}$. For the remaining partial waves, since the resonances dominate any change due to the polarization effect, the humps in the phase shift and corresponding small peaks in the delay are masked.

\subsection{Average time delay}
\label{Sec: Average time delay}
\begin{figure}
\includegraphics[width=12.5 cm]{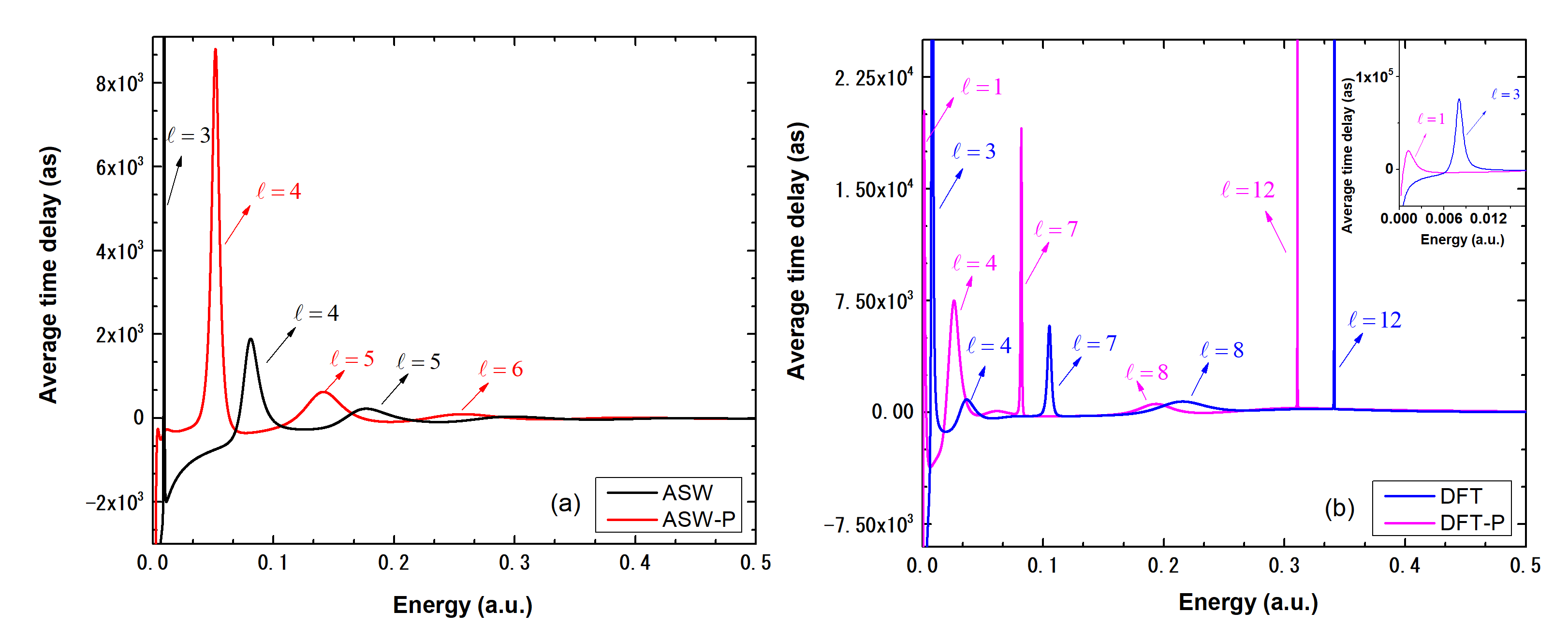}
\caption{ Average time delay of $e-C_{60}$ elastic scattering using model potentials (a) ASW (black) and ASW-P (red) (b) DFT (blue) and DFT-P (magenta). The average time delay behavior magnified in the low energy limit is shown in the inset for the DFT case. The resonant time delays are labeled with the corresponding partial wave's angular momentum quantum number. \label{fig8}}
\end{figure}   
The time delay results presented above are for the respective partial waves, which can not be experimentally observed. However, the average time delay is the quantity that can possibly be measured. Figure \ref{fig8} shows the average time delay computed using equation (\ref{avg_td}) in $e-C_{60}$ elastic scattering for (a) ASW, ASW-P and (b) DFT, DFT-P cases. The average time delay profiles exhibit all of the resonant peaks that were seen in TCS, as expected. For ASW and ASW-P cases (figure \ref{fig8}(a)), while approaching energy $E \to 0$, the average time delay tends to the negative infinity. The $s-$wave time delay profile dominates average time delay behavior in the low energy limit, since the low-energy collision is dominated by the s-wave scattering. For both the ASW and ASW-P cases, a progressive decline in the time delay peak value is noted. The contribution of partial wave $\ell=1$ is responsible for the very small hump in the ASW-P average time delay profile that is observed in the energy limit closer to zero. The average time delay in DFT also contains all of the resonant partial wave peaks (figure \ref{fig8}(b)), much like in the ASW case study. Because of the dominant influence of the $\ell=0$ partial wave, the average time delay value, here also, approaches negative infinity in the low energy limit, in accordance with the Wigner threshold law. We recall, the partial wave $\ell=12$ yielded the largest EWS time delay value in both the DFT and DFT-P cases. Accordingly, the average time delay also shows very steep peaks for the same partial wave. Structures from a common set of partial waves are observed for both DFT and DFT-P frames where the polarization effect causes an enhancement of the time delay and an energy dependent red-shift of the peak position is identified. 

\section{Summary and Conclusion}
\label{conculusion}

Calculation and analysis of the average and partial wave time delay for the low energy elastic $e-C_{60}$ scattering using ASW and DFT model potentials of fullerene are addressed in the present work. The impact of polarization effects on the time delay at these energies is also showcased. 

DFT potential, having a deeper well, shows a larger number of resonances compared to the ASW case. The Fano parameterization method is used to obtain the resonant energy values and other characteristic parameters. Since the polarization potential is attractive, the angular momentum barrier height gets lowered, resulting in a slightly lower resonant energy values for the polarization-added potentials. Corresponding to a resonance, a peak in the time delay is noted, which hints at the fact that the projectile-target duo stays for a longer time in their mutual vicinity at the resonance. This may be a useful feature while looking for the ideal targets for quantum memory applications \cite{lvovsky2009optical,ma2017optical}.

Minor peaks in the EWS time delay for resonant partial wave channels are observed as a consequence of the polarization effect, a result being reported here for the first time. Furthermore, the sensitivity of the time delay to changes in the interaction potential compared to other scattering parameters is evident from the present analysis. A similar observation was already made earlier for the case of photoionization, where the angular distribution asymmetry parameter ($\beta$) showed less sensitivity to certain interactions, whereas the time delay was extremely sensitive \cite{baral20226p, deshmukh2014attosecond}. The present work reaffirms this observation but for the first time in the context of electron scattering. 

 The present calculations employ realistic $e-C_{60}$ interaction potential using a sophisticated DFT methodology in combination with a PWA model. The resonant time delays are generally found to be comparatively longer relative to the background scattering. Although the DFT $e-C_{60}$ interaction is sophisticated compared to the static ASW model, there are realistic molecular effects that are missing in the present investigation. However, because of the symmetry of the $C_{60}$, it has been proven beyond doubt that static model potentials are simple but robust methods to simulate the $e-C_{60}$ interaction \cite{winstead2006elastic}. Moreover, at such lower impact-energies, the far longer projectile wavelength will not effectively "see" the target's atomistic details.   
 
 Using ultrafast two photon pump-probe laser pulses, measurements of time delay in photoionization are now possible \cite{schultze2010delay}. If experimental technology permits, it may also be possible to employ similar approaches to measure the time delay in electron scattering. An electron pulse can be sent to the target to scatter in a spherical pulse (scattering pulse), which can be followed by a probe laser pulse time-delayed from the original pulse in a controllable manner. This general technique can be realized either on an interferometric or a streaking track, as routinely done nowadays in the photoelectron chronoscopy. The current study may generate a preliminary impetus to this goal. \\

\textbf{Acknowledgements:}
The research is supported by DST-SERB-CRG Project No: \\
CRG/2022/000191, India (J.J.) and  CRG/2022/002309, India (H.R.V), and by the US
National Science Foundation Grant No. PHY-2110318 (H.S.C.).


\bibliography{reference}

\end{document}